\documentclass{amsart}
\usepackage{setspace}
\usepackage{a4}
\usepackage{amssymb,amsmath,amsthm,latexsym}
\usepackage{amsfonts}
\usepackage{amsfonts}
\usepackage{graphicx}
\usepackage{textcomp}
\usepackage{cite}
\usepackage{enumerate}
\usepackage[mathscr]{euscript}
\usepackage{mathtools}
\newtheorem{theorem}{Theorem}[section]

\newtheorem{definition}[theorem]{Definition}

\newtheorem{problem}[theorem]{Problem}
\newtheorem{proposition}[theorem]{Proposition}
\newtheorem{remark}[theorem]{Remark}

\title{This is the title}
\usepackage{amssymb}
\usepackage{amssymb}
\usepackage{amssymb}
\usepackage{amssymb}
\usepackage{amsmath}
\usepackage{tikz}
\usepackage{hyperref}
\usepackage{enumerate}
\usepackage{mathtools}
\usepackage{amsmath}
\usepackage{tikz}
\usepackage{amssymb}
\usepackage{amsmath}
\usepackage{tikz}
\usepackage{hyperref}
\usepackage{amssymb}
\usepackage{amssymb}
\usepackage{amssymb}
\usepackage{amssymb}
\usepackage{amsmath}
\usepackage{tikz}
\usepackage{hyperref}
\usepackage{enumerate}
\usepackage{mathtools}
\usepackage{amsmath}
\usepackage{tikz}
\usepackage{graphicx}
\usepackage{textcomp}
\usetikzlibrary{cd}
\usepackage{fancyhdr}

\begin{document}
\begin{center}
{\bf{PAULSEN AND PROJECTION PROBLEMS   FOR BANACH SPACES}}\\
\textbf{K. MAHESH KRISHNA} \\
Post Doctoral Fellow \\
Statistics and Mathematics Unit\\
Indian Statistical Institute, Bangalore Centre\\
Karnataka 560 059 India\\
Email: kmaheshak@gmail.com \\
Date: \today
\end{center}

\hrule
\vspace{0.5cm}
\textbf{Abstract}: Based on  the two decades old celebrated Paulsen problem and its solutions for Hilbert spaces by Kwok, Lau, Lee,  Ramachandran, Hamilton, and Moitra, we formulate Paulsen problem for Banach spaces. We also formulate projection problem for Banach spaces.

\textbf{Keywords}:  Paulsen problem, Projection problem, Approximate Schauder Frame.

\textbf{Mathematics Subject Classification (2020)}: 42C15.\\

\hrule
\tableofcontents

\section{Introduction}
Infinite dimensional frame theory, originated from the work of Duffin and Schaeffer \cite{DUFFINSCHAEFFER, OLE}  influenced the development of frame theory for finite dimensional Hilbert spaces towards the end of 20th century and theory developed rapidly in the first decade of this century. Revolutionary works on finite frames are \cite{CASAZZAFICKUSAP, CASAZZACLA, CASAZZACUSTOM, CAHILLCASAZZA2345, BALANCASAZZAEDIDIN, CASAZZAKOVACEVIC1234, BENEDETTOFICKUS, CASAZZAEDIDINEQUIVALENTS, DYKEMA, BODMANNPAULSEN, STRAWN, DYKEMAFREEMAN, STROHMERHEATH, CASAZZALEONHARD}. Finite dimensional definition of frames reads as follows. Letter $\mathcal{H}$ always denotes a finite dimensional Hilbert space. We denote the identity operator on $\mathcal{H}$ by $I_\mathcal{H}.$
\begin{definition} \cite{HANKORNELSONLARSONWEBER, BENEDETTOFICKUS}
	A collection $\{\tau_j\}_{j=1}^n$ in  a finite dimensional Hilbert space $\mathcal{H}$ over $\mathbb{K} $ ($\mathbb{R} $ or $\mathbb{C} $) is said to be a \textbf{frame} for 	$\mathcal{H}$ if there exist  $a,b>0$ such that 
	\begin{align*}
		a\|h\|^2\leq \sum_{j=1}^{n}|\langle h,\tau_j \rangle |^2 \leq b \|h\|^2, \quad \forall h \in \mathcal{H}.
	\end{align*}
	A frame $\{\tau_j\}_{j=1}^n$ for   $\mathcal{H}$ is said to be \textbf{Parseval} if 
	\begin{align*}
		h= \sum_{j=1}^{n}\langle h,\tau_j \rangle \tau_j , \quad \forall h \in \mathcal{H}.
	\end{align*}
A frame $\{\tau_j\}_{j=1}^n$ for   $\mathcal{H}$ is said to be \textbf{tight}  if 
\begin{align*}
	a\|h\|^2= \sum_{j=1}^{n}|\langle h,\tau_j \rangle |^2, \quad \forall h \in \mathcal{H}.
\end{align*}
A frame $\{\tau_j\}_{j=1}^n$ for   $\mathcal{H}$ is said to be \textbf{unit norm}   if $\|\tau_j\|=1$, for all $1\leq j\leq n$.
\end{definition}
A frame $\{\tau_j\}_{j=1}^n$ for $\mathcal{H}$ induces three operators defined as follows. 
\begin{enumerate}[\upshape (i)]
	\item Analysis operator $\theta_\tau: \mathcal{H} \ni h \mapsto \theta_\tau h \coloneqq (\langle h, \tau_j \rangle ) _{j=1}^n\in \mathbb{K}^n$.
	\item Synthesis operator $\theta_\tau^*: \mathbb{K}^n \ni (a_j)_{j=1}^n \mapsto \theta_\tau^*(a_j)_{j=1}^n\coloneqq\sum_{j=1}^{n}a_j\tau_j \in \mathcal{H}$.
	\item Frame operator $S_\tau: \mathcal{H} \ni h \mapsto S_\tau h \coloneqq \sum_{j=1}^{n}\langle h,\tau_j \rangle \tau_j \in \mathcal{H}$.
\end{enumerate}
Frame condition then gives that the analysis operator is injective, synthesis operator is surjective and the frame operator is positive and invertible. Further, the frame operator factors as $S_\tau=\theta_\tau^*\theta_\tau$. Since frame operator is positive and invertible,  one has the conclusion that $\{S_{\tau}^{-1/2}\tau_j\}_{j=1}^n$  is a Parseval frame for $\mathcal{H}$. In general, given a frame, construction of  Parseval frame in this way is difficult. Also, given a collection it is easily to verify that whether it is a frame than to verify that it is a Parseval frame. This leads to the following notions.  
\begin{definition}\cite{CAHILLCASAZZA} \label{ENF}
	A Parseval frame $\{\tau_j\}_{j=1}^n$ for a $d$-dimensional Hilbert space $\mathcal{H}$  is called an  \textbf{equal norm Parseval frame} if 
	\begin{align*}
		\|\tau_j\|^2=\frac{d}{n}, \quad \forall 1\leq j \leq n.
	\end{align*}
\end{definition}
\begin{definition}\cite{CAHILLCASAZZA}\label{NP}
	A  frame $\{\tau_j\}_{j=1}^n$ for a $d$-dimensional Hilbert space $\mathcal{H}$  is called an  \textbf{$\varepsilon$-nearly Parseval frame} ($\varepsilon<1$) if 	
	\begin{align*}
		(1-\varepsilon)I_\mathcal{H} \leq S_\tau \leq (1+\varepsilon)I_\mathcal{H}. 
	\end{align*}
\end{definition}
\begin{definition}\cite{CAHILLCASAZZA} \label{EENF}
	A  frame $\{\tau_j\}_{j=1}^n$ for a $d$-dimensional Hilbert space $\mathcal{H}$  is called an  \textbf{$\varepsilon$-nearly equal norm  frame} ($\varepsilon<1$) if 
	\begin{align*}
		(1-\varepsilon)\frac{d}{n}\leq 	\|\tau_j\|^2\leq (1+\varepsilon)\frac{d}{n}, \quad \forall 1\leq j \leq n.
	\end{align*}	
\end{definition}
\begin{definition}\cite{CAHILLCASAZZA}\label{BOTH}
	A frame for a $d$-dimensional Hilbert space $\mathcal{H}$ which is both $\varepsilon$-nearly  equal norm and  $\varepsilon$-nearly  Parseval is called as  \textbf{$\varepsilon$-nearly  equal norm Parseval frame.}
\end{definition}
Distance between two collections in Hilbert space are measured using following distance.
\begin{definition}\cite{CAHILLCASAZZA}\label{DIS}
	\textbf{Distance} between  two collections $\{\tau_j\}_{j=1}^n$,  $\{\omega_j\}_{j=1}^n$ in a Hilbert space $\mathcal{H}$ is defined as 
	\begin{align*}
		\operatorname{dist}(\{\tau_j\}_{j=1}^n, \{\omega_j\}_{j=1}^n)\coloneqq \left(\sum_{j=1}^{n}\|\tau_j-\omega_j\|^2\right)^\frac{1}{2}.
	\end{align*}
\end{definition}
We then have the fundamental Paulsen problem which is originated from the work of Holmes and Paulsen \cite{HOLMESPAULSEN}. Basic motivation is the following three results.
\begin{theorem}\cite{HOLMESPAULSEN}\label{PAULSENALGORITHM}
	There is an algorithm for turning every frame into an equal norm frame with same frame operator. 
\end{theorem}
\begin{theorem}\cite{CASAZZACUSTOM}\label{CLOSESTEQUALNORMFRAME}
 Let $\{\tau_j\}_{j=1}^n$ be an $\varepsilon$-nearly equal norm frame for  a $d$-dimensional Hilbert space $\mathcal{H}$. Then the closest equal norm frame to $\{\tau_j\}_{j=1}^n$ is given by $\{\omega_j\}_{j=1}^n$, where	
 \begin{align*}
 	\omega_j\coloneqq \left(\frac{\sum_{k=1}^{n}\|\tau_k\|}{n}\right)\frac{\tau_j}{\|\tau_j\|}, \quad \forall 1\leq j \leq n.
 \end{align*}
\end{theorem}
\begin{theorem}\cite{BODMANNCASAZZA, CASAZZAWERTY}\label{CLOSESTPARSEVALFRAME}
Let $\{\tau_j\}_{j=1}^n$ be a frame for  a $d$-dimensional Hilbert space $\mathcal{H}$. Then $\{S_\tau^{-1/2}\tau_j\}_{j=1}^n$ is the closest Parseval frame to $\{\tau_j\}_{j=1}^n$, i.e., 
\begin{align*}
	\sum_{j=1}^{n}\left\|S_\tau^\frac{-1}{2}\tau_j-\tau_j\right\|^2=\inf\left\{ \sum_{j=1}^{n}\left\|\tau_j-\omega_j\right\|^2: \{\omega_j\}_{j=1}^n \text{ is a Parseval frame for }	\mathcal{H}\right \}.
\end{align*}
Further, $\{S_\tau^{-1/2}\tau_j\}_{j=1}^n$ is the unique minimizer. Moreover, if $\{\tau_j\}_{j=1}^n$   is any $\varepsilon$-nearly  Parseval frame for $\mathcal{H}$, then 
\begin{align*}
		\sum_{j=1}^{n}\left\|S_\tau^\frac{-1}{2}\tau_j-\tau_j\right\|^2\leq d(2-\varepsilon-2\sqrt{1-\varepsilon})\leq \frac{d\varepsilon^2}{4}.
\end{align*}
\end{theorem}
\begin{problem}\cite{CAHILLCASAZZA} (\textbf{Paulsen problem})\label{PAUSENFIRST}
\textbf{Find the function  $f: (0, 1)\times \mathbb{N} \times \mathbb{N} \to [0, \infty)$ so that for any 	$\varepsilon$-nearly  equal norm  Parseval frame $\{\tau_j\}_{j=1}^n$ for d-dimensional Hilbert space $\mathcal{H}$, there is an equal norm Parseval frame  $\{\omega_j\}_{j=1}^n$ for $\mathcal{H}$ satisfying 
	\begin{align*}
		\operatorname{dist}^2(\{\tau_j\}_{j=1}^n, \{\omega_j\}_{j=1}^n)= \sum_{j=1}^{n}\|\tau_j-\omega_j\|^2\leq f(\varepsilon, n,d).
	\end{align*}
	Moreover, whether $f$ depends on $n$?}
\end{problem}
It is clear that Problem  \ref{PAUSENFIRST} can also be stated as follows.
\begin{problem}\cite{CAHILLCASAZZA}(\textbf{Paulsen problem})\label{PAULSENSECOND}
	\textbf{Find the function $f: (0, 1)\times \mathbb{N} \times \mathbb{N} \to [0, \infty)$ so that for any  $\varepsilon$-nearly  equal norm Parseval frame $\{\tau_j\}_{j=1}^n$  for d-dimensional Hilbert space $\mathcal{H}$, 
	\begin{align*}
		\inf \{\operatorname{dist}^2(\{\tau_j\}_{j=1}^n, \{\omega_j\}_{j=1}^n):\{\omega_j\}_{j=1}^n \text{ is an equal norm Parseval frame for } \mathcal{H}\}\leq f(\varepsilon, n,d).
	\end{align*}
	Moreover, whether $f$ depends on $n$?}
\end{problem}
The function $f$ in Definition \ref{PAUSENFIRST} and in Definition \ref{PAULSENSECOND} is known as Paulsen function. One first asks whether such a function exists. Using compactness of the unit sphere in finite dimensional Hilbert space, Hadwin proved the following.
\begin{theorem}\cite{BODMANNCASAZZA}\label{BODMANNCASAZZAEXISTS}
	Solution $f$	to Paulsen problem exists.
\end{theorem}
Casazza \cite{CASAZZABOOK} noticed that Paulsen function is bounded below which is independent of the number of elements in the frame.
\begin{proposition}\cite{CASAZZABOOK}
	Paulsen function $f$ satisfies 
	\begin{align*}
		f(\varepsilon, n, d)\geq \varepsilon^2d, \quad \forall \varepsilon>0, \forall d \in \mathbb{N}.
	\end{align*}	
\end{proposition}
Using the system of ordinary differential equations and using the notion of frame energy, Bodmann and Casazza made the first breakthrough for Paulsen problem in 2010 \cite{BODMANNCASAZZA}.
\begin{theorem}\cite{BODMANNCASAZZA}
	Let $n$ and $d$ be relatively prime and let $\varepsilon<1/2$. If $\{\tau_j\}_{j=1}^n$  is any $\varepsilon$-nearly equal norm frame  for a $d$-dimensional Hilbert space $\mathcal{H}$, then there is an equal norm Parseval frame  $\{\omega_j\}_{j=1}^n$ for $\mathcal{H}$ such that 
	\begin{align*}
			\operatorname{dist}^2(\{\tau_j\}_{j=1}^n, \{\omega_j\}_{j=1}^n)\leq \frac{29}{8}d^2n(n-1)^8\varepsilon.
	\end{align*}
\end{theorem}

Using the notion of Hilbert-Schmidt norm of operator, Casazza, Fickus, and Mixon \cite{CASAZZAFICKUSMIXON} gave the following version of Paulsen problem.
\begin{problem}\cite{CASAZZAFICKUSMIXON} (\textbf{Paulsen problem})\label{CASAZZASOLUTION}
\textbf{Given $n,d \in \mathbb{N}$, find positive 	$\delta(n,d)$, $c(n,d)$ and $\alpha(n,d)$ such that given any unit norm frame $\{\tau_j\}_{j=1}^n$ for $d$-dimensional Hilbert space $\mathcal{H}$ satisfying 
		\begin{align*}
	\left\|S_\tau-\frac{n}{d}I_\mathcal{H}\right\|_{\operatorname{HS}}\leq \delta(n,d)	,
	\end{align*}
there exists a unit norm tight frame $\{\omega_j\}_{j=1}^n$ for $\mathcal{H}$ such that 
	\begin{align*}
	\|\theta_\omega^*-\theta_\tau^*\|_{\operatorname{HS}} \leq c(n,d)	\left\|S_\tau-\frac{n}{d}I_\mathcal{H}\right\|_{\operatorname{HS}}^{\alpha(n,d)}	,
	\end{align*}
where HS denotes the Hilbert-Schmidt norm of the operator}.
\end{problem}
Using gradient of the frame potential, following particular cases   of Problem  \ref{CASAZZASOLUTION} have been solved in 2012 \cite{CASAZZAFICKUSMIXON}.
\begin{theorem}\cite{CASAZZAFICKUSMIXON}\label{PRIME}
	Let $n$ and $d$ be relatively prime and let $0<t<1/{2n}$.	Let $\{\tau^{(0)}_j\}_{j=1}^n$ be a unit norm tight frame for a $d$-dimensional Hilbert space $\mathcal{H}$ such that 
		\begin{align*}
		\left\|S^{(0)}_\tau-\frac{n}{d}I_\mathcal{H}\right\|_{\operatorname{HS}}^2\leq \frac{2}{d^3}	.
	\end{align*}
Define $\{\omega_j^{(k)}\}_{j=1}^n$ as follows.
\begin{align*}
\omega_j^{(k)}\coloneqq S_\tau^{(k)}\tau_j^{(k)}-\langle S_\tau^{(k)}\tau_j^{(k)}, \tau_j^{(k)}\rangle \tau_j^{(k)}, \quad \forall 1 \leq j \leq n, \forall k \geq 0.	
\end{align*}
Now define $\{\tau_j^{(k)}\}_{j=1}^n$ as follows.
\begin{align*}
\tau_j^{(k+1)}\coloneqq 
\begin{dcases}
	\cos(\|\omega_j^{(k)}\|t)\tau_j^{(k)}- \sin(\|\omega_j^{(k)}\|t)\frac{\omega_j^{(k)}}{\|\omega_j^{(k)}\|}& \text{if }  \omega_j^{(k)}\neq 0 \\
\tau_j^{(k)} & \text{if }  \omega_j^{(k)}=0\\
\end{dcases}
	,\quad \forall 1 \leq j \leq n, \forall k \geq 0.
\end{align*}
Then the limit of $\{\tau_j^{(k)}\}_{j=1}^n$ as $k\to \infty$ exists, denoted by $\{\tau_j^{(\infty)}\}_{j=1}^n$ is a unit norm tight frame for $\mathcal{H}$ and satisfies
\begin{align*}
	\left\|S^{(\infty)}_\tau-S^{(0)}_\tau\right\|_{\operatorname{HS}}\leq \frac{4d^{20}n^{8.5}}{1-2nt}	\left\|S^{(0)}_\tau-\frac{n}{d}I_\mathcal{H}\right\|_{\operatorname{HS}}.
\end{align*}
\end{theorem}
\begin{theorem}\cite{CASAZZAFICKUSMIXON}
Let $\varepsilon\leq 1/{2n}$. If $\{\tau_j\}_{j=1}^n$  is any $\varepsilon$-orthogonally partitionable unit norm frame  for a $d$-dimensional Hilbert space $\mathcal{H}$, then there is an orthogonally partitionable  unit  norm  frame  $\{\omega_j\}_{j=1}^n$ for $\mathcal{H}$ such that 	
\begin{align*}
\|\theta_\omega^*-\theta_\tau^*\|_{\operatorname{HS}} \leq \sqrt{2n}	(\varepsilon d)^\frac{1}{3}.
\end{align*}
\end{theorem}
\begin{theorem}\cite{CASAZZAFICKUSMIXON}
Let $n$ and $d$ be not relatively prime. 	If $\{\tau_j\}_{j=1}^n$  is any  unit norm frame  for a $d$-dimensional Hilbert space $\mathcal{H}$, then there is a unit  norm  frame  $\{\omega_j\}_{j=1}^n$ for $\mathcal{H}$ which is either tight or is orthogonally partitionable with equal redundancies in each of the two partitioned subsets such that 
\begin{align*}
\|\theta_\omega^*-\theta_\tau^*\|_{\operatorname{HS}} \leq	3 d^\frac{6}{7}\sqrt{n}	\left\|S_\tau-\frac{n}{d}I_\mathcal{H}\right\|_{\operatorname{HS}}^\frac{1}{7}.
\end{align*}
\end{theorem}
Projection problem is another classical problem in Hilbert space theory which states as follows.
\begin{problem}\cite{CAHILLCASAZZA} (\textbf{Projection problem})
\textbf{	Let $\mathcal{H}$ be a d-dimensional Hilbert space with orthonormal basis $\{u_k\}_{k=1}^d$. Find the function $g: (0, 1)\times \mathbb{N} \times \mathbb{N} \to [0, \infty)$ satisfying the following: If $P: \mathcal{H} \to \mathcal{H}$ is an orthogonal projection  of rank $n$ satisfying 
	\begin{align*}
		(1-\varepsilon)\frac{n}{d}\leq 	\|Pu_k\|^2\leq (1+\varepsilon)\frac{n}{d}, \quad \forall 1\leq k \leq d,
	\end{align*}
	then there exists an orthogonal projection $Q: \mathcal{H} \to \mathcal{H}$ with 
	\begin{align*}
		\|Qu_k\|^2=\frac{n}{d},	\quad \forall 1\leq k \leq d,
	\end{align*}
	satisfying 
	\begin{align*}
		\sum_{k=1}^{d}\|Pu_k-Qu_k\|^2\leq g(\varepsilon, n,d).	
	\end{align*}}
\end{problem}
Using chordal distance between subspaces \cite{CONWAYHARDIN}, Cahill and Casazza \cite{CAHILLCASAZZA}  showed that Paulsen problem has a solution if and only if  projection problem has solution.
\begin{theorem}\cite{CAHILLCASAZZA}
If $f$ is the function for the Paulsen problem and $g$ is the function for the projection problem, then 
\begin{align*}
g(\varepsilon, n,d)\leq 4 f(\varepsilon, n,d)\leq 8f(\varepsilon, n,d)		, \quad \forall \varepsilon, n, d.
\end{align*}
\end{theorem}
Using Naimark complement    of frames (see \cite{CASAZZAFICKUS2, CASAZZAKUTYNIOK, CASAZZALYNCH, CZAJA})  Cahill and Casazza \cite{CAHILLCASAZZA}  proved the following result. 
\begin{theorem}\cite{CAHILLCASAZZA}\label{COMPLEMENTIMPLIES}
	Let $n>d$. If $f$ is the function for the Paulsen problem, then 
	\begin{align*}
	f(\varepsilon, n,d) \leq 8 f\left(\frac{d}{n-d}, n, n-d\right).
	\end{align*}
\end{theorem}
Theorem \ref{COMPLEMENTIMPLIES} then gives the following result.
\begin{theorem}\cite{CAHILLCASAZZA}\label{LESSTHAN}
To solve the Paulsen problem for $d$-dimensional Hilbert space $\mathcal{H}$, it suffices to solve it for  Parseval frames $\{\tau_j\}_{j=1}^n$ for $\mathcal{H}$ with $d\leq n\leq 2d$.
\end{theorem}


In 2017, using operator scaling algorithm and smoothed analysis, Kwok, Lau, Lee, and Ramachandran \cite{KWOKLAULEERAMACHANDRAN, KWOKLAULEERAMACHANDRAN1} resolved Paulsen problem by deriving the following result.
\begin{theorem}\cite{KWOKLAULEERAMACHANDRAN, KWOKLAULEERAMACHANDRAN1}
	\label{RAMA}	For any  $\varepsilon$-nearly  equal norm Parseval frame $\{\tau_j\}_{j=1}^n$  for $\mathbb{R}^d$, there 	is an equal norm Parseval frame  $\{\omega_j\}_{j=1}^n$ for $\mathbb{R}^d$ satisfying 
	\begin{align*}
		\operatorname{dist}^2(\{\tau_j\}_{j=1}^n, \{\omega_j\}_{j=1}^n)\leq O(\varepsilon d^\frac{13}{2}).
	\end{align*}
	In other words, $f$ does not depend upon $n$.
\end{theorem}
In 2018, using radial isotropic position, Hamilton and Moitra \cite{HAMILTONMOITRA, HAMILTONMOITRA2} gave another proof of Paulsen problem which improved Theorem \ref{RAMA}.
\begin{theorem}\cite{HAMILTONMOITRA, HAMILTONMOITRA2}
	For any  $\varepsilon$-nearly  equal norm Parseval frame $\{\tau_j\}_{j=1}^n$  for $\mathbb{R}^d$, there 	is an equal norm Parseval frame  $\{\omega_j\}_{j=1}^n$ for $\mathbb{R}^d$ satisfying 
	\begin{align*}
		\operatorname{dist}^2(\{\tau_j\}_{j=1}^n, \{\omega_j\}_{j=1}^n)\leq 20 \varepsilon d^2.
	\end{align*}
	In other words,  $f(\varepsilon, n,d)=20\varepsilon d^2$.
\end{theorem}

\section{Paulsen and projection problems for Banach spaces}
Throughout this paper, we consider only finite dimensional Banach spaces. Let $\mathcal{X}$ be a $d$-dimensional Banach space. Given an operator $T$ on $\mathcal{X}$, by $\sigma(T)$, we mean the spectrum of $T$.

Motivated from the foundational work of Casazza, Han, and Larson \cite{CASAZZAHANLARSON} and from the fundamental work of Daubechies and DeVore \cite{DAUBECHIESDEVORE},  Casazza, Dilworth, Odell, Schlumprecht, and Zsak \cite{CASAZZADILWORTHODELL} introduced the notion of approximate Schauder frames for Banach spaces which led to the notion of approximate Schauder frames. Later, in the course of defining frame potential for Banach space and from the notion of Auerbach basis  \cite{WOJTASZCZYK, HAJEKBOOK, LORENZ, DIESTEL}, Chavez-Domingues, Freeman, and Kornelson \cite{KORNELSON} introduced the notion of unit norm tight frames for Banach spaces.  
\begin{definition}\cite{FREEMANODELL, THOMAS, CASAZZAHANLARSON, KORNELSON, KRISHNAJOHNSON}\label{ASFDEF} 
	Let $\{\tau_j\}_{j=1}^n$ be a collection in a Banach space  $\mathcal{X}$ and 	$\{f_j\}_{j=1}^n$ be a collection in  $\mathcal{X}^*$ (dual of $\mathcal{X}$). The pair $ (\{f_j \}_{j=1}^n, \{\tau_j \}_{j=1}^n) $ is said to be an \textbf{approximate Schauder frame (ASF)} for $\mathcal{X}$ if the frame operator 
	\begin{align*}
	S_{f, \tau}:\mathcal{X}\ni x \mapsto S_{f, \tau}x\coloneqq \sum_{j=1}^n
	f_j(x)\tau_j \in
	\mathcal{X}	
	\end{align*}
 is  invertible. Let $\lambda \in \mathbb{K}\setminus \{0\}$. An ASF  $ (\{f_j \}_{j=1}^n, \{\tau_j \}_{j=1}^n) $ is said to be \textbf{$\lambda$-tight} if $S_{f, \tau}=\lambda I_\mathcal{X}$, where $I_\mathcal{X}$ is the identity operator on $\mathcal{X}$. If $S_{f, \tau}= I_\mathcal{X}$, then we say  $ (\{f_j \}_{j=1}^n, \{\tau_j \}_{j=1}^n) $ is a \textbf{Parseval ASF} or \textbf{framing} or \textbf{Schauder frame} for  $\mathcal{X}$.  A tight ASF  $ (\{f_j \}_{j=1}^n, \{\tau_j \}_{j=1}^n) $ for   $\mathcal{X}$ is called as a \textbf{finite unit norm tight frame} if 
 \begin{align*}
 	\|f_j\|=\|\tau_j\|=f_j(\tau_j)=1, \quad \forall 1\leq j \leq n.
 \end{align*}
\end{definition} 
We now formulate Banach space versions of Definitions \ref{ENF}, \ref{NP}, \ref{EENF}  \ref{BOTH} and \ref{DIS}.
\begin{definition}\label{BENF}
	A Parseval frame $ (\{f_j \}_{j=1}^n, \{\tau_j \}_{j=1}^n) $ for a $d$-dimensional Banach space $\mathcal{X}$  is called an  \textbf{equal norm Parseval ASF} if 
	\begin{align*}
		\|\tau_j\|^2=\|f_j\|^2=f_j(\tau_j)=\frac{d}{n},\quad \forall 1\leq j \leq n.
	\end{align*}
\end{definition}
\begin{definition}\label{BNP}
	An ASF $ (\{f_j \}_{j=1}^n, \{\tau_j \}_{j=1}^n) $ for a $d$-dimensional Banach space $\mathcal{X}$  is called an  \textbf{$\varepsilon$-nearly Parseval ASF} ($\varepsilon<1$) if 	
	\begin{align*}
		\sigma(S_{f,\tau})\subseteq (1-\varepsilon, 1+\varepsilon).
	\end{align*}
\end{definition}
\begin{definition}\label{BEENF}
	An ASF   $ (\{f_j \}_{j=1}^n, \{\tau_j \}_{j=1}^n) $ for a $d$-dimensional Banach space $\mathcal{X}$  is called an  \textbf{$\varepsilon$-nearly equal norm  ASF} ($\varepsilon<1$) if 
		\begin{align*}
		(1-\varepsilon)\frac{d}{n}\leq 	\|\tau_j\|^2=f_j(\tau_j)=\|f_j\|^2\leq (1+\varepsilon)\frac{d}{n}, \quad \forall 1\leq j \leq n.
	\end{align*}	
\end{definition}
\begin{definition}\label{BBOTH}
	An ASF  which is both $\varepsilon$-nearly  equal norm and  $\varepsilon$-nearly  Parseval is called as  \textbf{$\varepsilon$-nearly  equal norm Parseval ASF}.
\end{definition}
\begin{definition}\label{BDIS}
	Let $\{\tau_j\}_{j=1}^n$, $\{\omega_j\}_{j=1}^n$ be  collections in a $d$-dimensional Banach space  $\mathcal{X}$ and 	$\{f_j\}_{j=1}^n$, $\{g_j\}_{j=1}^n$ be  collections in  $\mathcal{X}^*$. 	\textbf{Distance} between  two tuples $ (\{f_j \}_{j=1}^n, \{\tau_j \}_{j=1}^n) $ and $ (\{g_j \}_{j=1}^n, \{\omega_j \}_{j=1}^n) $ is defined as 
	\begin{align}\label{DISTANCE}
		\operatorname{dist}((\{f_j \}_{j=1}^n, \{\tau_j \}_{j=1}^n) ,  (\{g_j \}_{j=1}^n, \{\omega_j \}_{j=1}^n))\coloneqq \left(\sum_{j=1}^{n}\frac{1}{2}(\|\tau_j-\omega_j\|^2+\|f_j-g_j\|^2)\right)^\frac{1}{2}.
	\end{align}
\end{definition}
\begin{remark}
	\begin{enumerate}[\upshape(i)]
		\item 	Rather defining distance as in Definition (\ref{DISTANCE}), we would have also been defined distance as 
		\begin{align*}
			\operatorname{dist}^*((\{f_j \}_{j=1}^n, \{\tau_j \}_{j=1}^n) ,  (\{g_j \}_{j=1}^n, \{\omega_j \}_{j=1}^n))\coloneqq	\frac{1}{2}\left(\left(\sum_{j=1}^{n}\|\tau_j-\omega_j\|^2\right)^\frac{1}{2}+\left(\sum_{j=1}^{n}\|f_j-g_j\|^2\right)^\frac{1}{2}\right).
		\end{align*}
		However, we then have 
		\begin{align*}
			\operatorname{dist}^*((\{f_j \}_{j=1}^n, \{\tau_j \}_{j=1}^n) ,  (\{g_j \}_{j=1}^n, \{\omega_j \}_{j=1}^n))&\leq 	\operatorname{dist}((\{f_j \}_{j=1}^n, \{\tau_j \}_{j=1}^n) ,  (\{g_j \}_{j=1}^n, \{\omega_j \}_{j=1}^n)), \\&\quad \forall (\{f_j \}_{j=1}^n, \{\tau_j \}_{j=1}^n) ,  (\{g_j \}_{j=1}^n, \{\omega_j \}_{j=1}^n).
		\end{align*}
		\item We can also try by replacing the  distance in Definition (\ref{DISTANCE}) by
		\begin{align*}
			\operatorname{dist}_p((\{f_j \}_{j=1}^n, \{\tau_j \}_{j=1}^n) ,  (\{g_j \}_{j=1}^n, \{\omega_j \}_{j=1}^n))\coloneqq \left(\sum_{j=1}^{n}\frac{1}{2}(\|\tau_j-\omega_j\|^p+\|f_j-g_j\|^p)\right)^\frac{1}{p}.	
		\end{align*}
		for some $p>0$.
	\end{enumerate}
\end{remark}
It is clear that for a (Hilbert space) frame $\{\tau_j\}_{j=1}^n$ for a Hilbert space $\mathcal{H}$, by defining $f_j:\mathcal{H}\ni h \mapsto \langle h, \tau_j \rangle \in \mathbb{K}$, for all $1\leq j \leq n$, Definitions  \ref{BENF}, \ref{BNP}, \ref{BEENF},  \ref{BBOTH} and \ref{BDIS} reduce to Definitions \ref{ENF}, \ref{NP}, \ref{EENF},  \ref{BOTH} and \ref{DIS}, respectively. 
In view of Theorems \ref{PAULSENALGORITHM},   \ref{CLOSESTEQUALNORMFRAME}   and \ref{CLOSESTPARSEVALFRAME} we ask following three problems.
\begin{problem}
\textbf{Whether 	there is an algorithm for turning every ASF into an equal norm ASF with same frame operator?}
\end{problem}
\begin{problem}
\textbf{What is the closest (in terms of distance given in Definition \ref{BDIS}) Parseval ASF to  a given ASF?}
\end{problem}
\begin{problem}
\textbf{What is the closest equal norm ASF (in terms of distance given in Definition \ref{BDIS}) to a given equal norm ASF?}
\end{problem}
Here we state Paulsen problem for Banach spaces.
\begin{problem}(\textbf{Paulsen problem for Banach spaces})\label{PAULSENBANACH}
\textbf{Find the function  $f: (0, 1)\times \mathbb{N} \times \mathbb{N} \to [0, \infty)$ so that for any 	$\varepsilon$-nearly  equal norm  Parseval ASF $ (\{f_j \}_{j=1}^n, \{\tau_j \}_{j=1}^n) $ for a d-dimensional Banach  space $\mathcal{X}$, there is an equal norm Parseval ASF  $ (\{g_j \}_{j=1}^n, \{\omega_j \}_{j=1}^n) $  for $\mathcal{X}$ satisfying 
	\begin{align*}
		\operatorname{dist}^2((\{f_j \}_{j=1}^n, \{\tau_j \}_{j=1}^n) ,  (\{g_j \}_{j=1}^n, \{\omega_j \}_{j=1}^n))= \sum_{j=1}^{n}\frac{1}{2}(\|\tau_j-\omega_j\|^2+\|f_j-g_j\|^2)\leq f(\varepsilon, n,d).
	\end{align*}
	Moreover, whether $f$ depends on $n$?}	
\end{problem}
We can reformulate Problem \ref{PAULSENBANACH} as follows. 
\begin{problem}(\textbf{Paulsen problem for Banach spaces})
       	\textbf{Find the function $f: (0, 1)\times \mathbb{N} \times \mathbb{N} \to [0, \infty)$ so that for any  $\varepsilon$-nearly  equal norm Parseval ASF $ (\{f_j \}_{j=1}^n, \{\tau_j \}_{j=1}^n) $  for a d-dimensional Banach  space $\mathcal{X}$, 
       	\begin{align*}
       		\inf \{&\operatorname{dist}^2((\{f_j \}_{j=1}^n, \{\tau_j \}_{j=1}^n) ,  (\{g_j \}_{j=1}^n, \{\omega_j \}_{j=1}^n)): (\{g_j \}_{j=1}^n, \{\omega_j \}_{j=1}^n)\text{ is an equal norm Parseval }\\
       			&\text{ASF for } \mathcal{X}\}\leq f(\varepsilon, n,d).
       	\end{align*}
       	Moreover, whether $f$ depends on $n$?}	
\end{problem}
In view of Theorem \ref{BODMANNCASAZZAEXISTS}   we first ask the following problem.
\begin{problem}
\textbf{Does there exists a solution for the Paulsen problem for Banach spaces (Problem \ref{PAULSENBANACH})?}
\end{problem}
Here we state projection problem for Banach spaces.
\begin{problem} (\textbf{Projection problem for Banach spaces})\label{PROBANACH}
	\textbf{Let $\mathcal{X}$ be a $d$-dimensional Banach  space with Auerbach  basis $(\{\zeta_k\}_{k=1}^d, \{u_k\}_{k=1}^d)$  ($\{\zeta_k\}_{k=1}^d$ is the collection of dual functionals). Find the function $g: (0, 1)\times \mathbb{N} \times \mathbb{N} \to [0, \infty)$ satisfying the following: If $P: \mathcal{X} \to \mathcal{X}$ is a projection (idempotent) of rank $n$ satisfying 
	\begin{align*}
		(1-\varepsilon)\frac{n}{d}\leq 	\|Pu_k\|^2=	\|\zeta_kP\|^2=|\zeta_k(Pu_k)|\leq (1+\varepsilon)\frac{n}{d}, \quad \forall 1\leq k \leq d,
	\end{align*}
	then there exists a  projection (idempotent) $Q: \mathcal{X} \to \mathcal{X}$ with 
	\begin{align*}
		\|Qu_k\|^2=	\|\zeta_kQ\|^2=|\zeta_k(Pu_k)|=\frac{n}{d},	\quad \forall 1\leq k \leq d,
	\end{align*}
	satisfying 
	\begin{align*}
		\sum_{k=1}^{d}\frac{1}{2}(\|Pu_k-Qu_k\|^2+\|\zeta_kP-\zeta_kQ\|^2)\leq g(\varepsilon, n,d).	
	\end{align*}}
\end{problem}
We already remarked that Paulsen and projection problems were connected using chordal distance. Even though we are unable to connect Pausen and projection problems for Banach spaces, we introduce the notion of chordal distance for subspaces of Banach spaces. Chordal distance between two subspaces of a finite dimensional Hilbert space is defined using angles between them. Since we do not have inner product we take alternative way using a characterization  to introduce the notion for Banach spaces.
\begin{definition}
\textbf{Let $\mathcal{X}$ be a $d$-dimensional Banach  space and $P, Q: \mathcal{X}\to \mathcal{X}$ be rank $m$ projections onto subspaces $\mathcal{Y}$ and $\mathcal{Z}$ of $\mathcal{X}$, respectively. We define the chordal distance  between $\mathcal{Y}$ and $\mathcal{Z}$ as 
\begin{align*}
	\operatorname{dist}_{\operatorname{chordal}}(\mathcal{Y}, \mathcal{Z})\coloneqq (m-\operatorname{Trace}(PQ))^\frac{1}{2}.
\end{align*}}
\end{definition}

\begin{problem}
\textbf{Whether there is a relation between Paulsen problem for Banach spaces (Problem \ref{PAULSENBANACH}) and projection problem for Banach spaces (Problem \ref{PROBANACH})?}
\end{problem}
We next note that dilation result for a class of ASFs has been derived in \cite{KRISHNAJOHNSONDIL}. For finite dimensional spaces  it states as follows. We refer  \cite{KRISHNA, KRISHNAJOHNSONDIL} for the notion of approximate Riesz bases  for  finite dimensional Banach spaces.
\begin{theorem}\label{KRISHNAJOHNSON}\cite{KRISHNAJOHNSONDIL}
Let $ (\{f_j \}_{j=1}^n, \{\tau_j \}_{j=1}^n) $ be an ASF for   $\mathcal{X}$. Then there exists a Banach space $\mathcal{Y}$ which contains $\mathcal{X}$ isometrically and an approximate Riesz basis  $ (\{g_j \}_{j=1}^n, \{\omega_j \}_{j=1}^n) $  for   $\mathcal{Y}$ such that 
\begin{align*}
	f_j=g_jP_{|\mathcal{X}}, ~ \tau_j=P\omega_j, \quad \forall 1\leq j \leq n,
\end{align*}
where $P:\mathcal{Y}\to \mathcal{X}$ is an onto projection. If $ (\{f_j \}_{j=1}^n, \{\tau_j \}_{j=1}^n) $ is a Parseval ASF  or Schauder frame  for   $\mathcal{X}$, then $ (\{g_j \}_{j=1}^n, \{\omega_j \}_{j=1}^n) $  can be taken as $g_j(\omega_k)=\delta_{j,k}, \forall 1\leq j, k \leq n$.
\end{theorem}
Using Theorem  \ref{KRISHNAJOHNSON} and based on Theorem \ref{PRIME} we can ask following problem.
\begin{problem}
	\textbf{Can Paulsen problem for Banach spaces (Problem \ref{PAULSENBANACH}) be solved for $n$ and $d$ which are relatively prime?}
\end{problem}
In view of Theorem \ref{LESSTHAN} we ask the following problem.
\begin{problem}
	\textbf{Can the Paulsen problem for Banach spaces (Problem \ref{PAULSENBANACH})  be reduced like Theorem  \ref{LESSTHAN} using Theorem \ref{KRISHNAJOHNSON}?}
\end{problem}
We can also formulate the following problems.
\begin{problem}
	\textbf{Classify  Banach spaces for which the solution to  Paulsen problem for Banach spaces \ref{PAULSENBANACH}	exists.}
\end{problem}
\begin{problem}
		\textbf{Classify Banach spaces for which solution to projection  problem for Banach spaces \ref{PROBANACH} exists.}
\end{problem}
 \begin{problem}
 	\textbf{Classify Banach spaces for which both Paulsen problem \ref{PAULSENBANACH} and projection problem \ref{PROBANACH} have solutions.}
 \end{problem} 
  \begin{remark}
  	\textbf{A web  of conjectures and problems for p-approximate Schauder frames \cite{MAHESHTHESIS, KRISHNAJOHNSON} for Banach spaces have been made by the author  in \cite{KRISHNA}. Note that for finite collections in finite dimensional Banach spaces, the notions approximate Schauder frames and p-approximate Schauder frames coincide.}
  \end{remark}

 \bibliographystyle{plain}
 \bibliography{reference.bib}

\end{document}